\documentclass[leqno,a4paper,12pt]{article}

\usepackage{epsfig,amssymb,amsmath,amsthm}

\newcommand{\ds}{\displaystyle}

\numberwithin{equation}{section}
\newtheorem{theorem}{Theorem}[section]

\newtheorem{conjecture}{Conjecture}
\newtheorem{proposition}{Proposition}[section]

\begin{document}

\title{Families of (1,2)--Symplectic Metrics on Full 
Flag Manifolds}

\author{\sc{Marlio Paredes}\thanks{Supported by CAPES--Brazil and
COLCIENCIAS--Colombia}}

\date{}

\maketitle

\begin{abstract}
We obtain new families of (1,2)--symplectic invariant metrics on the
full complex flag manifolds $F(n)$. For $n \geq 5$, we characterize $n-3$
different $n$--dimensional families of (1,2)--symplectic invariant metrics on
$F(n)$. Any of these families corresponds to a different class of
non--integrable invariant
 almost complex structure on $F(n)$. 
 
\end{abstract}

\section{Introduction} \label{s:int}

Recently Mo and Negreiros \cite{MN1}, by using moving frames and tournaments,  
showed explicitly the existence of a $n$--dimensional family of invariant
(1,2)--symplectic metrics on $F(n) = \frac{U(n)}{U(1) \times \cdots \times
U(1)}$.  This family corresponds to the family of the parabolic almost
complex structures on $F(n)$. In this paper we study the existence of other
families of invariant (1,2)--symplectic metrics corresponding to classes of
non--integrable invariant  almost complex structure on $F(n)$ different to the
parabolic. 

Eells and Sampson \cite{ES}, proved that if $\phi \colon M \rightarrow
N$ is a holomorphic map between K\"ahler manifolds then $\phi$ is
harmonic. This result was generalized by Lichnerowicz (see \cite{L} or
\cite{Sa}) as follows: Let $(M,g,J_1)$ and $(N,h,J_2)$ be almost Hermitian 
manifolds with $M$ cosymplectic and $N$ (1,2)--symplectic. Then any
holomorphic map $\phi \colon (M,g,J_1) \rightarrow (N,h,J_2)$ is harmonic.

If we like to obtain harmonic maps, $\phi \colon M^2 \rightarrow
F(n)$, from a closed Riemannian surface $M^2$ to a full flag manifold 
$F(n)$, by the Lichnerowicz theorem we have to study (1,2)--symplectic
metrics on $F(n)$ because a Riemann surface is a K\"ahler
manifold and we know that a K\"ahler manifold is a cosymplectic manifold (see
\cite{Sa} or \cite{GH}).

To study the invariant Hermitian geometry  of $F(n)$ it is natural to begin by
studing its invariant almost complex structures. Borel and Hirzebruch \cite{BH}, 
proved that there are $2^{\binom{n}{2}}$ $U(n)$--invariant almost complex  
structures on $F(n)$. This number is the same number of tournaments
with $n$ players or nodes. A tournament is a digraph in which any two
nodes are joined by exactly one oriented edge (see \cite{M} or
\cite{BS}). There is a natural identification between almost complex  
structures on $F(n)$ and tournaments with $n$ players (see \cite{MN2} 
or \cite{BS}).

Tournaments can be classified in isomorphism classes. In this
classification, one of these classes corresponds to the integrable
structures and the other ones correspond to non--integrable
structures. Burstall and Salamon \cite{BS}, proved that a almost
complex structure $J$ on $F(n)$ is integrable if and only if the
associated tournament to $J$ is isomorphic to the canonical tournament
(the canonical tournament with $n$ players, $\{1,2,\ldots,n\}$, is
defined by $i \rightarrow j$ if and only if $i < j$).

Borel proved that exits a $(n-1)$--dimensional family
 of invariant K\"ahler
metrics on $F(n)$ for each invariant complex
 structure on $F(n)$ (see
\cite{Be} or \cite{Bo}). Eells and Salamon \cite{ESa}, proved that any
parabolic structure on $F(n)$ admits a (1,2)--symplectic metric. Mo and
Negreiros \cite{MN1}, showed explicitly that there is a $n$--dimensional 
family of invariant (1,2)--symplectic metrics for each parabolic
 structure on
$F(n)$. 
 
In this paper, we characterize new $n$--parametrical families of
(1,2)--sym\-plec\-tic invariant metrics on $F(n)$, different to the K\"ahler 
and parabolic. More precisely, we obtain explicitly $n-3$ different
$n$--dimensional families of (1,2)--symplectic invariant metrics, for
each $n \geq 5$. Each of them corresponds to a different class of
non--integrable
invariant almost complex structure on $F(n)$. These metrics
are used to produce new examples of harmonic maps $\phi
 \colon M^2
\rightarrow F(n)$, using the result of Lichnerowicz
above.
 
This paper is part of the author's Doctoral Thesis \cite{P1}. I wish 
to thank my advisor Professor Caio Negreiros for his right advise. I 
would like to thank Professor Xiaohuan
 Mo for his helpful comments and
dicussions on this work.

\section{Preliminaries}

A full flag manifold is defined by
\begin{equation} \label{eq:2.1}
F(n) = \{ (L_1, \dots, L_n) : L_i \ \text{is a subspace of} \
{\mathbb C}^n , \mbox{dim}_{\mathbb C}L_i = 1, \ \ L_i \bot L_j \}.
\end{equation}
The unitary group $U(n)$ acts transitively on $F(n)$. Using this
action we obtain an algebraic description for $F(n)$:
\begin{equation} \label{eq:2.2}
F(n) = \frac{U(n)}{T} = \frac{U(n)}{\underbrace{U(1) \times \cdots  
\times U(1)}_{n-times}} \ ,
\end{equation}
where $T = \underbrace{U(1) \times \cdots \times U(1)}_{n-times}$ is
a maximal torus in $U(n)$.

Let $\mathfrak{p}$ be the tangent space of $F(n)$ in $(T)$. An invariant 
almost complex structure on $F(n)$ is a linear map $J \colon \mathfrak{p}
\rightarrow \mathfrak{p}$ such that $J^2 = - I$.

A tournament or $n$--tournament ${\cal T}$, consists of a finite set 
$T= \{p_1, p_2,\dots, \\ p_n\}$ of $n$ players, together with a dominance 
relation, $\to$, that assigns to every pair of players a winner, 
i.e. $p_i \to p_j$ or $p_j \to p_i$. If $p_i \to p_j$ then we say 
that $p_i$ beats $p_j$. A tournament ${\cal T}$ may be represented by a
directed graph in which $T$ is the set of vertices and any two vertices are
joined by an oriented edge. 

Let ${\cal T}_1$ be a tournament with $n$ players $\{1, \ldots , n\}$
and ${\cal T}_2$ another tournament with $m$ players $\{1, \ldots ,
m\}$. A homomorphism between ${\cal T}_1$ and ${\cal T}_2$ is a
mapping $\phi: \{1, \dots , n\} \to \{1, \dots , m\}$ such that
\begin{equation} \label{eq:2.3}
s \overset{ {\cal T}_1 }{\longrightarrow} t \quad \Longrightarrow 
\quad \phi(s) \overset{ {\cal T}_2 }{\longrightarrow} \phi(t) 
\qquad \mbox{or} \qquad \phi(s) = \phi(t).
\end{equation}
When $\phi$ is bijective we said that ${\cal T}_1$ and ${\cal T}_2$
are isomorphic.

An $n$--tournament determines a score vector $(s_1, \dots , s_n)$,
such that $\ds \sum_{i=1}^{n} s_i \\ = \binom{n}{2}$, with components
equal the number of games won by each player. Isomorphic tournaments have
identical score vectors. Figure  \ref{fig:1} shows the isomorphism classes of
$n$--tournaments for  $n = 2,3,4$, together with their score vectors. This
Figure was taken of Moon's book \cite{M}. In the Moon's notation not all of
the arcs are included in the drawings, if an arc joining two nodes has not
been drawn then it is to be understood that the arc is oriented from the
higher node to the lower node.

\begin{figure}
\begin{center}
\input{fig1famil.pstex_t}
\end{center}
\caption{Isomorphism classes of $n$--tornaments to $n = 2,3,4$.}
\label{fig:1}
\end{figure}

The canonical $n$--tournament ${\cal T}_n$ is defined by setting 
$i \rightarrow j$ if and only if $i < j$. Up to isomorphism, 
${\cal T}_n$ is the unique $n$--tournament  satisfying the following
equivalent conditions:
\begin{itemize}
\item the dominance relation is transitive, i.e. if $i \to j$ \ and \
$j \to k$ \ then \ $i \to k$,
\item there are no 3--cycles, i.e. closed paths $i_1 \to i_2 \to \to i_3 
\to i_1$, see \cite{M},
\item the score vector is $(0,1,2,\ldots,n-1)$.
\end{itemize}

For each invariant almost complex structure $J$ on $F(n)$, we can 
associate a $n$--tour\-na\-ment ${\cal T}(J)$ in the following way:
If $J(a_{ij}) = (a'_{ij})$ then ${\cal T}(J)$ is such that for $i < j$ 
\begin{equation} \label{eq:2.4}
\Bigl( i \to j \ \Leftrightarrow \ a'_{ij} = \sqrt{-1} \, a_{ij} \Bigl) 
\quad \text{or} \quad \Bigl( i \leftarrow j \ \Leftrightarrow  
\ a'_{ij} = - \sqrt{-1} \, a_{ij} \Bigl),
\end{equation}
see \cite{MN2}.

An almost complex structure $J$ on $F(n)$ is said to be integrable if 
$F(n)$ is a complex manifold, i.e. $F(n)$ admits complex coordinate
systems with holomorphic coordinate changes. Burstall and Salamon
\cite{BS} proved the following result:
\begin{theorem} \label{t:2.1}
An almost complex structure $J$ on $F(n)$ is integrable if and only
if ${\cal T}(J)$ is isomorphic to the canonical tournament ${\cal T}_n$.
\end{theorem}
Thus, if ${\cal T}(J)$ contains a 3--cycle then $J$ is not integrable.
Classes (2) and (4) in Figure \ref{fig:1} correspond
to the integrable almost
complex structures on $F(3)$ and $F(4)$ respectively.

An invariant almost complex structure $J$ on $F(n)$ is called
parabolic if there is a permutation $\tau$ of $n$ elements such that the
associate tournament ${\cal T}(J)$ is given, for $i < j$, by
$$
\Bigl( \tau(j) \to \tau(i), \quad \text{if} \ j-i \ \text{is even} \Bigl)
\qquad \text{or} \qquad \Bigl( \tau(i) \to \tau(j), \quad \text{if} \ j-i \ 
\text{is odd} \Bigl).
$$
Classes (3) and (7) in Figure \ref{fig:1} represent the parabolic structures 
on $F(3)$ and  $F(4)$ respectively.
    
A $n$--tournament ${\cal T}$, for $n \geq 3$, is called irreducible or
Hamiltonian if it contains a $n$--cycle, i.e. a path $\pi(n) \to \pi(1) \to
\pi(2) \to \cdots \to \pi(n-1) \to \pi(n)$, where $\pi$ is a permutation of
$n$ elements.

A $n$--tournament ${\cal T}$ is transitive if given three nodes $i,j,k$ of 
${\cal T}$ then $i \to j \ \text{and} \ j \to k \
\Longrightarrow \ i \to k$. The canonical tournament is the only one
transitive tournament up to isomorphisms.

We consider $\mathbb{C}^n$ equipped with the standard Hermitian 
inner product, i.e. for $V = (v_1,\ldots,v_n)$ and 
$W = (w_1,\ldots,w_n)$ \ in \ $\mathbb{C}^n$, we have $\left< V,W \right> =
\sum_{i=1}^n v_i \overline{w_i}$. We use the convention $\overline{v_i} =
v_{\bar\imath} \ \mbox{and} \ \overline{f_{i\bar{\jmath}}} =
f_{\bar{\imath}j}$.

A frame consists of an ordered set of $n$ vectors $(Z_1,\ldots,Z_n)$,
such that $Z_1 \wedge \ldots \wedge Z_n \neq 0$, and it is called 
unitary, if $\left< Z_i , Z_j \right> = \delta_{i\bar{\jmath}}$. The 
set of unitary frames can be identified with the unitary group. 

If we write $dZ_i = \sum_j \omega_{i\bar{\jmath}} Z_j$, the coefficients
$\omega_{i\bar{\jmath}}$ are the Maurer--Cartan forms of the unitary group
$U(n)$. They are skew--hermitian, i.e. $\omega_{i\bar{\jmath}} + 
\omega_{\bar{\jmath}i} = 0$. For more details see \cite{ChW}.

We may define all left invariant metrics on $(F(n),J)$ by (see \cite{Bl} 
or \cite{N1})
\begin{equation} \label{eq:2.5}
ds^2_{\Lambda} = \sum_{i,j} {\lambda}_{ij} \omega_{i\bar{\jmath}} 
\otimes \omega_{\bar{\imath}j},
\end{equation} 
where $\Lambda = ({\lambda}_{ij})$ is a real matrix such that:
\begin{equation} \label{eq:2.6}
{\lambda}_{ij} \left \{ 
\begin{array}{lcl}
> 0, & \mbox{se} & i \neq j \\
= 0, & \mbox{se} & i = j
\end{array}
\right., 
\end{equation}
and the Maurer--Cartan forms $\omega_{i\bar{\jmath}}$ are such that
\begin{equation} \label{eq:2.7}
\omega_{i\bar{\jmath}} \in {\mathbb C}^{1,0} \
\mbox{((1,0) type forms)} \quad \Longleftrightarrow  \quad 
i \overset{{\cal T}(J)}{\longrightarrow} j .
\end{equation}

The metrics (\ref{eq:2.5}) are called Borel type and they are almost
Hermitian for every invariant almost complex structure $J$, i.e. 
$ds^2_{\Lambda} (JX , JY) = ds^2_{\Lambda} (X , Y)$, for all tangent 
vectors $X,Y$. When $J$ is integrable $ds^2_{\Lambda}$ is said to be
Hermitian.

Let $J$ be an invariant almost complex structure on $F(n)$,
${\cal T}(J)$ the associated tournament, and $ds^2_{\Lambda}$ an
invariant metric. The K\"ahler form with respect to $J$ and 
$ds^2_{\Lambda}$ is defined by 
\begin{equation} \label{eq:2.8}
\Omega(X,Y) = ds^2_{\Lambda} (X , JY),
\end{equation}
for any tangent vectors $X,Y$. For each permutation $\tau$, of $n$ elements, 
the K\"ahler form can be
write in the following way (see \cite{MN1})
\begin{equation} \label{eq:2.9}
\Omega = -2\sqrt{-1} \ \sum_{i<j} \mu_{\tau (i) \tau (j)} 
\omega_{\tau (i) \overline{\tau (j)}} \wedge 
\omega_{\overline{\tau (i)}\tau (j)},
\end{equation}
where $\mu_{\tau (i) \tau (j)} = \varepsilon_{\tau (i) \tau (j)} 
\lambda_{\tau (i) \tau (j)}$ \quad and \quad $ \varepsilon_{ij} = \left \{ 
\begin{array}{rcl}
1  & \mbox{se} & i \to j \\
-1 & \mbox{se} & j \to i \\
0  & \mbox{se} & i = j
\end{array}
\right.$.

Let $J$ be an invariant almost complex structure on $F(n)$. Then
$F(n)$ is said to be almost K\"ahler if and only if $\Omega$ is
closed, i.e. $d\Omega = 0$. If $J$ is integrable and $\Omega$ is
closed then $F(n)$ is said to be a K\"ahler manifold.

Mo and Negreiros proved in \cite{MN1}
\begin{equation} \label{eq:2.10}
d\Omega = 4 \sum_{i<j<k} C_{\tau (i) \tau (j) \tau (k)} 
\Psi_{\tau (i) \tau (j) \tau (k)}, 
\end{equation}
where
\begin{equation} \label{eq:2.11}
C_{ijk} = \mu_{ij} - \mu_{ik} + \mu_{jk}, 
\end{equation}
and
\begin{equation} \label{eq:2.12}
\Psi_{ijk} = \mbox{Im} (\omega_{i \bar{\jmath}} \wedge 
\omega_{\bar{\imath}k} \wedge \omega_{j\bar{k}}).
\end{equation}

We denote by $\mathbb{C}^{p,q}$ the space of complex forms with degree
$(p,q)$ on $F(n)$. Then, for any $i,j,k$, we have either
$\Psi_{ijk} \in \mathbb{C}^{0,3} \oplus \mathbb{C}^{3,0} \quad \text{or} 
\quad \Psi_{ijk} \in \mathbb{C}^{1,2} \oplus \mathbb{C}^{2,1}$.

An invariant almost Hermitian metric $ds^2_{\Lambda}$ is said to be
(1,2)--symplectic if and only if $(d\Omega)^{1,2} = 0$. If
$d^*\Omega=0$ then the metric is said to be cosymplectic.

The following result due to Mo and Negreiros \cite{MN1}, is very useful
to study (1,2)--symplectic metrics on $F(n)$:

\begin{theorem} \label{t:2.2}
If $J$ is a $U(n)$--invariant almost complex structure on $F(n)$, 
$n \geq 4$, such that ${\cal T}(J)$ contains one of 4--tournaments (5) or (6)
in the Figure \ref{fig:1} then $J$ does not admit any invariant
(1,2)--symplectic metric.
\end{theorem}

\section{Main Theorem}

It is known that, on $F(3)$ there is a 2--parametric family of K\"ahler
metrics and a 3--parametric family of (1,2)--symplectic metrics
corresponding to the non--integrable almost complex structures class (the
parabolic class). Then each invariant almost complex structure on $F(3)$
admits a (1,2)--symplectic metric, see \cite{ESa}, \cite{Bo}.

On $F(4)$ there are four isomorphism classes of 4--tournaments or
equivalently almost complex structures and the Theorem \ref{t:2.2}
shows that two of them do not admit (1,2)--symplectic metric. The another
two classes correspond to the K\"ahler and parabolic cases. $F(4)$
has a 3--parametric family of K\"ahler metrics and a 4--parametric
family of (1,2)--symplectic metrics which is not K\"ahler, see
\cite{MN1}.

On $F(5)$, \ $F(6)$ \ and \ $F(7)$ we have the following families of
(1,2)--symplectic invariant metrics, different to the K\"ahler and parabolic:
on $F(5)$, two 5--parame\-tric families; on $F(6)$, four 6--parametric
families, two of them generalizing the two families on $F(5)$ and, on $F(7)$
there are eight 7--para\-me\-tric families, four of them generalizing the four
ones on $F(6)$, this results are contained in \cite{P1} or \cite{P2}. 

\begin{figure}
\begin{center}
\input{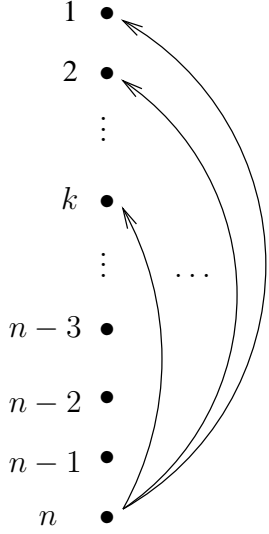}
\caption{Tournaments of the Theorem \ref{t:3.1}}
\label{fig:2}
\end{center}
\end{figure}

In this paper we proved the following result:

\begin{theorem} \label{t:3.1}
Let $J$ be an invariant almost complex structure on $F(n)$ such 
that the associated tournament ${\cal T}(J)$ is the
tournament in Figure \ref{fig:2}. An invariant metric
$ds^2_{\Lambda}$ is (1,2)-symplectic with respect to $J$ if and only
if the matrix $\Lambda = (\lambda_{ij})$ satisfies
$$
\lambda_{ij} = \lambda_{i(i+1)} + \lambda_{(i+1)(i+2)} + \cdots + 
\lambda_{(j-1)j}
$$
for $i = 1, \dots , n-1$ \ and \ $j = 2, \dots , n$, except to 
$\lambda_{1n}, \lambda_{2n}, \ldots , \lambda_{kn}$
which satisfy the following relations
$$
\begin{array}{rcl}
\lambda_{2n} & = & \lambda_{12} + \lambda_{1n}  \\
\lambda_{3n} & = & \lambda_{12} + + \lambda_{23} + \lambda_{1n} \\
             & \vdots &  \\
\lambda_{kn} & = &  \lambda_{12} + \lambda_{23} + \ldots + 
            \lambda_{(k-1)k} +  \lambda_{1n}.
\end{array}
$$
\end{theorem}

This theorem provides a $n$--family of (1,2)-- symplectic metrics on $F(n)$,
for each $1 \leq k \leq n-3$. This families are different to the family described 
by Mo and Negreiros in \cite{MN1} and corresponding to non--integrable almost
complex structures. All of our families are $n$-parametric. 

None of our families contains the normal metric. This fact is according
to the Wolf and Gray results in \cite{WG} which proved that the normal
metric on $F(n)$ is (1,2)-symplectic if and only if $n \leq 3$.

The score vector of these families can be write as:
$$
(1, \ 2, \ldots, \ k, \ k, \ldots, \ n-k-1, \ n-k-1, \ldots, 
\ n-3, \ n-2) \, ,
$$
for $n \geq 2k + 1$.

In order to prove this theorem we present in the following  section some
preliminary results.

\section{The Families for $k$ = 1, \, 2, \, 3, \, 4 }

\begin{figure}
\begin{center}
\input{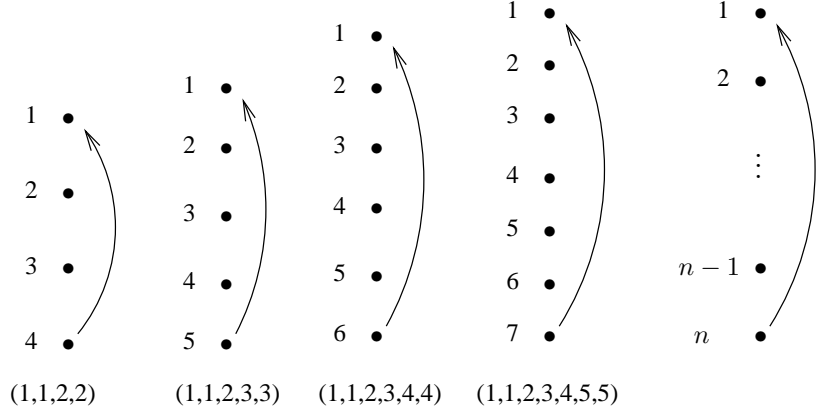}
\caption{Tournaments of the family for $k=1$}
\label{fig:3}
\end{center}
\end{figure}

\begin{proposition} \label{p:4.1}
Let $J$ be an invariant almost complex structure on $F(n)$, with 
$n \geq 4$, such that the associated tournament ${\cal T}(J)$ is the
last tournament in the Figure \ref{fig:3}. An invariant metric
$ds^2_{\Lambda}$ is (1,2)-symplectic with respect to $J$ if and only
if the matrix $\Lambda = (\lambda_{ij})$ satisfies
$$
\lambda_{ik} = \lambda_{i(i+1)} + \lambda_{(i+1)(i+2)} + \cdots + 
\lambda_{(k-1)k}
$$
for $i = 1, \dots , n-1$ \ and \ $k = 2, \dots , n$, except to 
$\lambda_{1n}$.
\end{proposition}

The corresponding matrix ${\Lambda}^1$ has the way:
{\scriptsize
$$
\Lambda^1 =
\left(
\begin{array}{cccccc}
0 & \lambda_{12} & \lambda_{12} + \lambda_{23} & \cdots & \lambda_{12} + 
\cdots + \lambda_{(n-2)(n-1)} & \lambda_{1n} \\

\lambda_{12} & 0 & \lambda_{23} & \dots & \lambda_{23} + \cdots + 
\lambda_{(n-2)(n-1)} &  \lambda_{23} + \cdots + \lambda_{(n-1)n}  \\

\lambda_{12} + \lambda_{23} & \lambda_{23} & 0 & \ddots & \vdots & 
\vdots \\

  * &   *   &  *  &  \cdots  & \lambda_{(n-2)(n-1)} & 
          \lambda_{(n-2)(n-1)} + \lambda_{(n-1)n} \\

 *   &   *  &  *  & \cdots &   0   & \lambda_{(n-1)n} \\
 *   &   *  &  *  & \cdots &  \lambda_{(n-1)n} &  0   
\end{array}
\right).
$$}

\begin{proof}
The proof is made by using induction over $n$. First we prove the result for 
$n = 4$, in this case the tournament ${\cal T}(J)$ is isomorphic to the
first tournament in Figure \ref{fig:3}. Calculating $d\Omega$ using
(\ref{eq:2.10}) we obtain
$$
\begin{array}{rcl}
d\Omega & = & C_{123} \Psi_{123} + C_{124} \Psi_{124} + C_{134} \Psi_{134} + 
              C_{234} \Psi_{234} \\
        & = & (\lambda_{12} - \lambda_{13} + \lambda_{23}) \Psi_{123} + 
              (\lambda_{12} + \lambda_{14} + \lambda_{24}) \Psi_{124} + \\
        &   & + (\lambda_{13} + \lambda_{14} + \lambda_{34}) \Psi_{134} +
              (\lambda_{23} - \lambda_{24} + \lambda_{34}) \Psi_{234} 
\end{array}
$$
and $d\Omega^{(1,2)} = (\lambda_{12} - \lambda_{13} + \lambda_{23}) \Psi_{123}
+  (\lambda_{23} - \lambda_{24} + \lambda_{34}) \Psi_{234}$. Then
$ds^2_{\lambda}$  is (1,2)--symplectic if and only if
$$
\left\{
\begin{array}{rcl}
\lambda_{12} - \lambda_{13} + \lambda_{23} & = & 0 \\
\lambda_{23} - \lambda_{24} + \lambda_{34} & = & 0
\end{array} \right\} 
\Longleftrightarrow 
\left\{
\begin{array}{rcl}
\lambda_{13} & = & \lambda_{12} + \lambda_{23} \\
\lambda_{24} & = & \lambda_{23} + \lambda_{34}
\end{array} \right\}.
$$

Suppose that the result is true to $n-1$. For $n$ we must consider 
two cases:

\begin{itemize}
\item[(a)] $i<j<k, \quad i \neq 1 \quad \text{or} \quad k \neq n$.
Then $\varepsilon_{ij} = \varepsilon_{ik} = \varepsilon_{jk} = 1$ \
and \ $C_{ijk} = \lambda_{ij} - \lambda_{ik} + \lambda_{jk} \neq 0$.
 
\item[(b)] $1<j<n$. Then $\varepsilon_{1j} = \varepsilon_{jn} = 1$, 
$\varepsilon_{1n} = -1$ \ and \ $C_{1jn} = \lambda_{1j} + \lambda_{1n} +
\lambda_{jn} \neq 0$. 
\end{itemize}
$$
\begin{array}{rcl}
\text{(a)} & \Rightarrow & (d\Omega)^{2,1} + (d\Omega)^{1,2} = 
\ds\sum_{i<j<k} C_{ijk} \Psi_{ijk}, \quad i \neq 1, \quad  k \neq n. \\
           &            &                    \\
\text{(b)} & \Rightarrow & (d\Omega)^{3,0} + (d\Omega)^{0,3} = 
\ds\sum_{j=2}^{n-1} C_{1jn} \Psi_{1jn} \neq 0.
\end{array}
$$
Then $ds^2_{\Lambda}$ is (1,2)-symplectic if and only if
$\Lambda = (\lambda_{ij})$ satisfies the linear system

{\tiny
$$
\begin{array}{rcl}
         \lambda_{12} - \lambda_{13} + \lambda_{23}  &   =    & 0 \\
         \lambda_{12} - \lambda_{14} + \lambda_{24}  &   =    & 0 \\
                                                     & \vdots &   \\
\lambda_{12} - \lambda_{1(n-1)} + \lambda_{2(n-1)}   &   =    & 0 \\
         \lambda_{13} - \lambda_{14} + \lambda_{34}  &   =    & 0 \\
                                                     & \vdots &   \\
\lambda_{13} - \lambda_{1(n-1)} + \lambda_{3(n-1)}   &   =    & 0 \\
         \lambda_{14} - \lambda_{15} + \lambda_{45}  &   =    & 0 \\
                                                     & \vdots &   \\
\lambda_{1(n-2)} - \lambda_{1(n-1)} + \lambda_{(n-2)(n-1)}    & = & 0 \\
         \lambda_{23} - \lambda_{24} +\lambda_{34}   &   =    & 0 \\
                                                     & \vdots &   \\
         \lambda_{23} - \lambda_{2n} +\lambda_{3n}   &   =    & 0 \\
                                                     & \vdots &   \\
 \lambda_{(n-3)(n-2)} - \lambda_{(n-3)n} + \lambda_{(n-2)n}    & = & 0  \\
 \lambda_{(n-2)(n-1)} - \lambda_{(n-2)n} + \lambda_{(n-1)n}    & = & 0. 
\end{array}
$$}

This system contains all of equations corresponding to the system for
$n-1$. Then all of elements of $\Lambda^1$ to $n-1$ are equal to the
matrix for $n$, except $\lambda_{1(n-1)}$.  Using the system above we
see how to write $\lambda_{1(n-1)}, \lambda_{2n}, \lambda_{3n}, \dots , 
\lambda_{(n-2)n}$:
{\footnotesize
$$
\begin{array}{rcl}
\lambda_{12} - \lambda_{1(n-1)} + \lambda_{2(n-1)} = 0 & \Longrightarrow 
& \lambda_{1(n-1)} = \lambda_{12} + \lambda_{2(n-1)}  \\
                  & \Longrightarrow &
\lambda_{1(n-1)} = \lambda_{12} + \lambda_{23} + \cdots + 
\lambda_{(n-2)(n-1)} \\
\lambda_{(n-2)(n-1)} - \lambda_{(n-2)n} + \lambda_{(n-1)n} = 0 & 
\Longrightarrow & \lambda_{(n-2)n} = \lambda_{(n-2)(n-1)} + 
\lambda_{(n-1)n} \\
        &                 &   \\
\lambda_{(n-3)(n-2)} - \lambda_{(n-3)n} + \lambda_{(n-2)n} = 0 & 
\Longrightarrow & \lambda_{(n-3)n} = \lambda_{(n-3)(n-2)} + 
\lambda_{(n-2)n} \\
          & \Longrightarrow & \lambda_{(n-3)n} = \lambda_{(n-3)(n-2)} +
            \lambda_{(n-2)(n-1)} + \mbox{} \\
          &  & \phantom{\lambda_{(n-3)n} =} +  \lambda_{(n-1)n} \\
          & \vdots &      \\
\lambda_{23} - \lambda_{2n} +\lambda_{3n} = 0 & \Longrightarrow & 
\lambda_{2n} = \lambda_{23} + \lambda_{3n} \\
       & \Longrightarrow & \lambda_{2n} = \lambda_{23} + 
         \lambda_{34} + \cdots + \lambda_{(n-1)n}.
\end{array}
$$}
\end{proof}

In $F(4)$ this family is the same as the family obtained by Mo and  
Negreiros \cite{MN1}, because the corresponding 4-tournament is
the parabolic 4-tournament. Any tournament of this family is irreducible
and such that any 4-subtournament of it is transitive (class (4) in Figure
\ref{fig:1}) or irreducible (class (7) in Figure \ref{fig:1}).

The following propositions are presented without proof. These propositions 
are proved in a similar way like was proved the proposition \ref{p:4.1}.

\begin{figure}
\begin{center}
\input{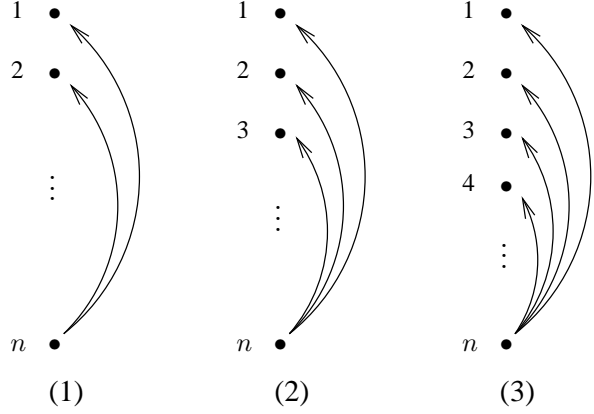}
\caption{Tournaments of the propositions \ref{p:4.2},  \ref{p:4.3} and
\ref{p:4.4}} \label{fig:4}
\end{center}
\end{figure}

\begin{proposition} \label{p:4.2}
Let $J$ be an invariant almost complex structure on $F(n)$, for 
$n \geq 5$, such that the associated tournament ${\cal T}(J)$ is the
tournament (1) in the Figure \ref{fig:4}. An invariant metric
$ds^2_{\Lambda}$ is (1,2)-symplectic with respect to $J$ if and only
if the matrix $\Lambda = (\lambda_{ij})$ satisfies
$$
\lambda_{ik} = \lambda_{i(i+1)} + \lambda_{(i+1)(i+2)} + \cdots + 
\lambda_{(k-1)k}
$$
for $i = 1, \dots , n-1$ \ and \ $k = 2, \dots , n$, except to 
$\lambda_{1n}$ and $\lambda_{2n}$ which satisfy $\lambda_{2n} = 
\lambda_{12} + \lambda_{1n}$.
\end{proposition}

In this case, the corresponding matrix $\Lambda^2$ is
{\tiny
$$
\Lambda^2 =
\left(
\begin{array}{cccccc}
0 & \lambda_{12} & \lambda_{12} + \lambda_{23} & \cdots & \lambda_{12} 
+ \cdots + \lambda_{(n-2)(n-1)} & \lambda_{1n} \\

\lambda_{12} & 0 & \lambda_{23} & \cdots & \lambda_{23} + \cdots + 
\lambda_{(n-2)(n-1)} &  \lambda_{12} + \lambda_{1n}  \\

\lambda_{12} + \lambda_{23} & \lambda_{23} & 0 & \ddots & \vdots & 
\vdots \\

  * &   *   &  *  &  \ddots  & \lambda_{(n-2)(n-1)} & 
  \lambda_{(n-2)(n-1)} + \lambda_{(n-1)n} \\

 *   &   *  &  *  & \cdots &   0   & \lambda_{(n-1)n} \\
 *   &   *  &  *  & \cdots &  \lambda_{(n-1)n} &  0   
\end{array}
\right).
$$}

\begin{proposition} \label{p:4.3}
Let $J$ be an invariant almost complex structure on $F(n)$, for 
$n \geq 6$, such that the associated tournament ${\cal T}(J)$ is the
tournament (2) in the Figure \ref{fig:4}. An invariant metric
$ds^2_{\Lambda}$ is (1,2)-symplectic with respect to $J$ if and only
if the matrix $\Lambda = (\lambda_{ij})$ satisfies
$$
\lambda_{ik} = \lambda_{i(i+1)} + \lambda_{(i+1)(i+2)} + \cdots + 
\lambda_{(k-1)k}
$$
for $i = 1, \dots , n-1$ \ and \ $k = 2, \dots , n$, except to 
$\lambda_{1n}$, $\lambda_{2n}$ and $\lambda_{3n}$ which satisfy 
$\lambda_{2n} = \lambda_{12} + \lambda_{1n}$ and $\lambda_{3n} = 
\lambda_{12} + \lambda_{23} + \lambda_{1n}$.
\end{proposition}

\begin{proposition} \label{p:4.4}
Let $J$ be an invariant almost complex structure on $F(n)$, for 
$n \geq 7$, such that the associated tournament ${\cal T}(J)$ is the
tournament (3) in the Figure \ref{fig:4}. An invariant metric
$ds^2_{\Lambda}$ is (1,2)-symplectic with respect to $J$ if and only
if the matrix $\Lambda = (\lambda_{ij})$ satisfies
$$
\lambda_{ik} = \lambda_{i(i+1)} + \lambda_{(i+1)(i+2)} + \cdots + 
\lambda_{(k-1)k}
$$
for $i = 1, \dots , n-1$ \ and \ $k = 2, \dots , n$, except to 
$\lambda_{1n}$, $\lambda_{2n}$, $\lambda_{3n}$ and $\lambda_{4n}$ 
which satisfy $\lambda_{2n} = \lambda_{12} + \lambda_{1n}$, \ 
$\lambda_{3n} = \lambda_{12} + \lambda_{23} + \lambda_{1n}$, and 
$\lambda_{4n} = \lambda_{12} + \lambda_{23} + \lambda_{34} + 
\lambda_{1n}$
\end{proposition}

 The corresponding matrices $\Lambda^3$ and $\Lambda^4$ is on the
final pages of this  paper. Any tournament of these families is irreducible
and such that any 4-subtournament of it is transitive (class (4) in Figure
\ref{fig:1}) or irreducible (class (7) in Figure \ref{fig:1}).

\section{Proof of the Main Theorem}

We use induction over $n$, begining with $n=4$. The proposition \ref{p:4.1}
shows that the result is true for $n = 4$. Suppose that the result is true for
$n-1$.

We need to calculate the coefificients $C_{ijk}$ in the formula
(\ref{eq:2.10}), then we have three types of 3-subtournaments of 
${\cal T}(J)$ to consider:

\begin{itemize}

\item[(a)]  To the 3-cycles we have
$$
C_{ijn} = \lambda_{ij} + \lambda_{in} + \lambda_{jn} \neq 0
$$
for $k<j<n$ \ and \ $i = 1, \ldots , k$. It implies that
$(d\Omega)^{3,0} \neq 0$.

\item[(b)] To the 3-subtournaments, $(ijn)$, such that $i< j \leq k$ \ and \ 
$i = 1, 2, \ldots , \\ k-1$, we have
$$
C_{ijn} = \lambda_{ij} + \lambda_{in} - \lambda_{jn}.
$$

\item[(c)] To the 3-subtournaments which do not satisfy neither (a) nor (b),
we have 
$$
C_{ijl} = \lambda_{ij} - \lambda_{il} + \lambda_{jl}, \quad i<j<l.
$$
\end{itemize}

(b) and (c) give us the information to calculate $(d\Omega)^{1,2}$.
Then the metric $ds^2_{\Lambda}$ is (1,2)-symplectic if and only if
the matrix $\Lambda = (\lambda_{ij})$ satisfies 

$$
\left\{
\begin{array}{lcl}
(d) \quad \lambda_{ij} + \lambda_{in} - \lambda_{jn} = 0; 
         & \quad & i< j \leq k, \quad i = 1, 2, \ldots , k-1 \\
         &  &  \\
(e) \quad \lambda_{ij} - \lambda_{il} + \lambda_{jl} = 0; 
         & \quad & i<j<l, \quad \mbox{do not satisfy (a) and (b)}.
\end{array}
\right.
$$

(d) and (e) include all of equations corresponding to the case for
$n-1$, except the equations given by the following 3-subtournaments
$$
(ij(n-1)), \qquad \mbox{with} \quad i = 1, \ldots , k-1, \quad j = 2, 
\ldots k, \quad \mbox{and} \quad i<j.
$$

Therefore, by the hipothesis induction, all of the elements of the
matrix $\Lambda^k$ corresponding to $n-1$ are equal to the matrix for
$n$, except the elements $\lambda_{1(n-1)}, \ \lambda_{2(n-1)}, \ldots , 
\lambda_{k(n-1)}$. Then we must calculate $\lambda_{1(n-1)}, \ldots , 
\lambda_{k(n-1)}, \\ \lambda_{2n}, \ldots,\lambda_{(n-2)n}$.

\begin{itemize}

\item[(i)] We take $i=k$, \ $j=k+1$ \ and \ $l=n-1$ in (e). Then
$$
\lambda_{k(k+1)} - \lambda_{k(n-1)} + \lambda_{(k+1)(n-1)} = 0,
$$
and it implies
$$
\begin{array}{rcl}
\lambda_{k(n-1)} & = & \lambda_{k(k+1)} + \lambda_{(k+1)(n-1)} \\
                 & = & \lambda_{k(k+1)} + \lambda_{(k+1)(k+2)} + 
		 \ldots + \lambda_{(n-2)(n-1)}.
\end{array}
$$
Using (e) again, with $i=k-1$, \ $j=k$ \ and \ $l=n-1$, we obtain
$$
\lambda_{(k-1)k} - \lambda_{(k-1)(n-1)} + \lambda_{k(n-1)} = 0,
$$
and it implies
$$
\begin{array}{rcl}
\lambda_{(k-1)(n-1)} & = & \lambda_{(k-1)k} + \lambda_{k(n-1)} \\
                     & = & \lambda_{(k-1)k} + \lambda_{k(k+1)} + \ldots +
                       \lambda_{(n-2)(n-1)}.
\end{array}
$$
If we continue using (e) to the rest of values: $i = k-2, \ldots, 2,1$, 
\ $j = k-1, \ldots, 2,1$ and $l=n-1$; we arrive at the following
equations  
$$
\begin{array}{c}
\lambda_{23} - \lambda_{2(n-1)} + \lambda_{3(n-1)} = 0 \\
\lambda_{12} - \lambda_{1(n-1)} + \lambda_{2(n-1)} = 0, 
\end{array}
$$
which imply
$$
\begin{array}{rcl}
\lambda_{2(n-1)}  & = &  \lambda_{23} + \lambda_{3(n-1)} \\
                  & = & \lambda_{23} + \lambda_{34} + \ldots + 
                        \lambda_{(n-2)(n-1)} \\
\end{array}
$$
and
$$
\begin{array}{rcl}                 
\lambda_{1(n-1)}  & = &  \lambda_{12} + \lambda_{2(n-1)} \\
                  & = & \lambda_{12} + \lambda_{23} + \ldots + 
                        \lambda_{(n-2)(n-1)}.
\end{array}
$$
Hence, the equation (e) implies
$$
\lambda_{i(n-1)} = \lambda_{i(i+1)} + \lambda_{(i+1)(i+2)} + \ldots +
                   \lambda_{(n-2)(n-1)}
$$
for $i = 1,2, \ldots , k$.

\item[(ii)] If $i=1$ and $j=2$ in (d) then $\lambda_{12} + 
\lambda_{1n} - \lambda_{2n} = 0$, and $\lambda_{2n} = \lambda_{12} + 
\lambda_{1n}$. Using again (d) with $i=1$ and $j=3$ we obtain 
$\lambda_{3n} = \lambda_{12} + \lambda_{23} + \lambda_{1n}$. We use
(d) repeatedly up to obtain
$$
\lambda_{in} = \lambda_{12} + \lambda_{23} + \ldots + \lambda_{(i-1)i} +
               \lambda_{1n}
$$
for $i= 2,3, \ldots , k$.

\item[(iii)] In order to calculate $\lambda_{(k+1)n}, \ldots , 
\lambda_{(n-2)n}$, we use (e) with $i = k+1,\ldots, \\ n-2$. We
obtain
{\tiny
$$
\begin{array}{rcrcl}
\lambda_{(n-2)(n-1)} - \lambda_{(n-2)n} + \lambda_{(n-1)n} = 0 & 
\Longrightarrow & \lambda_{(n-2)n} & = & \lambda_{(n-2)(n-1)} + 
\lambda_{(n-1)n} \\
    &     &     &    &   \\
\lambda_{(n-3)(n-2)} - \lambda_{(n-3)n} + \lambda_{(n-2)n} = 0 & 
\Longrightarrow & \lambda_{(n-3)n} & = & \lambda_{(n-3)(n-2)} + 
\lambda_{(n-2)n} \\
                       & \Longrightarrow &
\lambda_{(n-3)n} & = & \lambda_{(n-3)(n-2)} + \lambda_{(n-2)(n-1)} + \\ 
     &   &  &  &  \mbox{} + \lambda_{(n-1)n} \\
                             & \vdots &   &    &   \\
\lambda_{(k+1)(k+2)} - \lambda_{(k+1)n} +\lambda_{(k+2)n} = 0 
   & \Longrightarrow & \lambda_{(k+1)n} & = & \lambda_{(k+1)(k+2)} + 
                                          \lambda_{(k+2)n} \\
   & \Longrightarrow & \lambda_{(k+1)n} & = & \lambda_{(k+1)(k+2)} + 
                 \lambda_{(k+2)(k+3)} + \\
   &  &  &  & \mbox{} + \cdots + \lambda_{(n-1)n}.
\end{array}
$$} \hfill $\square$
\end{itemize}

\section{Harmonic Maps}

In this section we construct new examples of harmonic maps using the 
following result due to Lichnerowicz \cite{L}:

\begin{theorem} \label{t:6.1}
Let $\phi \colon (M,g) \rightarrow (N,h)$ be a $\pm$ holomorphic map between 
almost Hermitian manifolds where $M$ is cosymplectic and $N$ is
(1,2)-symplec\-tic. Then $\phi$ is harmonic.
\end{theorem}

In order to construct harmonic maps $\phi \colon M^2 \to F(n)$ using the
theorem above, we need to know examples of holomorphic maps. Then we use the
following construction due to Eells and Wood \cite{EW}.

Let $h: M^2 \to \mathbb{CP}^{n-1}$ be a full holomorphic map ($h$ is
full if $h(M)$ is not contained in none $\mathbb{CP}^k$, for all 
$k < n-1$). We can lift $h$ to $\mathbb{C}^n$, i.e. for every $p \in
M$ we can find a neighborhood of $p$, $U \subset M$, such that  
$h_U = (u_0, \dots , u_{n-1}): M^2 \supset U \to \mathbb{C}^n - {0}$
satisfies $h(z) = [h_U(z)] = [(u_0(z), \dots , u_{n-1}(z))]$.

We define the $k$-th associate curve of $h$ by
$$
\begin{array}{rcll}
{\cal O}_k: & M^2 & \longrightarrow & \mathbb{G}_{k+1}(\mathbb{C}^n)  \\
            &  z  & \longmapsto     & h_U(z) \wedge \partial h_U (z) \wedge
                     \dots \wedge \partial^k h_U (z),
\end{array}
$$
for $0 \leq k \leq n-1$. And we consider
$$
\begin{array}{rcll}
h_k: & M^2 & \longrightarrow & \mathbb{CP}^{n-1} \\
     &  z  & \longmapsto     & {\cal O}_k^{\perp} (z) 
                               \cap {\cal O}_{k+1} (z),
\end{array}
$$
for $0 \geq k \geq n-1$.

The following theorem, due to Eells and Wood \cite{EW}, is very important 
because it gives the classification of the harmonic maps from $S^2 \thicksim 
\mathbb{CP}^1$ into a projective space $\mathbb{CP}^{n-1}$.

\begin{theorem} \label{t:6.2}
For each $k \in \mathbb{N}, \ 0 \leq k \leq n-1$, $h_k$ is harmonic.
Furthermore, given $\phi: (\mathbb{CP}^1 , \ g) \to 
(\mathbb{CP}^{n-1} , \ Killing \ metric)$ a full harmonic map, then there
are unique $k$ and $h$ such that $\phi = h_k$.  
\end{theorem}

This theorem provides in a natural way the following holomorphic maps
$$
\begin{array}{rcll}
\Psi: & M^2 & \longrightarrow & F(n) \\
      &  z  & \longmapsto     & (h_0(z), \dots , h_{n-1}(z)),
\end{array}
$$
called by Eells-Wood's map (see \cite{N2}).

We can write the set of (1,2)-symplectic metrics on $F(n)$, 
characterized in the sections above, in the following way
$$
\mathfrak{M}_n = \{ g^k= ds^2_{\Lambda^k} : 1 \leq k \leq n-3 \}.
$$

Using Theorem \ref{t:6.1} we obtain the following results

\begin{proposition} \label{p:6.1}
Let $\phi: M^2 \to \bigl(F(n), g \bigl)$, $g \in \mathfrak{M}_n$ a
holomorphic map. Then $\phi$ is harmonic. 
\end{proposition} 

\begin{proposition} \label{p:6.2}
Let $\phi: (F(l), g) \to (F(n), \Tilde{g})$ a holomorphic map with 
$g \in \mathfrak{M}_l$ \ and \ $\Tilde{g} \in \mathfrak{M}_n$. Then
$\phi$ is harmonic.
\end{proposition} 

\section{Conjectures and Problems}

We would like to obtain some classification of the (1,2)-symplectic metrics
on $F(n)$. Our results suggest the following conjectures:

\begin{conjecture}
Let $(F(n),ds^2_{\Lambda},J)$ be a full flag manifold, with $n \geq 5$. The
Borel metric $ds^2_{\Lambda}$ is (1,2)-symplectic if and only if all of the
4-subtournament of associated tournament ${\cal T}(J)$ is transitive (class
(4) in Figure \ref{fig:1}) or irreducible (class (7) in Figure \ref{fig:1}).
\end{conjecture}

\begin{conjecture}
All of families of (1,2)-symplectic metrics on $F(n)$, different to
the K\"ahler metrics family, are n-dimensional.
\end{conjecture}

Another interesting problems are: 

\begin{itemize}
\item How many families of (1,2)-symplectic metrics on $F(n)$ are there?

\item Do exist co-symplectic metrics different to (1,2)-symplectic, on $F(n)$?

\end{itemize}

\newpage

\begin{center}
\rotatebox[origin=c]{90}{
\scriptsize
$
\Lambda^3 \ = \
\left(
\renewcommand{\arraystretch}{2.5}
\begin{array}{ccccccc}
0 & \lambda_{12} & \lambda_{12} + \lambda_{23} & \lambda_{12} + \lambda_{23} + \lambda_{34} & \cdots & \lambda_{12} + \cdots + \lambda_{(n-2)(n-1)} & \lambda_{1n} \\
%%%%%%%%%%%%%%%%%%%%%%%%%%%%%%%%%%%%%%%%%%%%%%%%%%%%%%%%%%%%%%%%%%%%%%%%%%%%%%
\lambda_{12} & 0 & \lambda_{23} & \lambda_{23} + \lambda_{34} & \cdots & \lambda_{23} + \cdots + \lambda_{(n-2)(n-1)} &  \lambda_{12} + \lambda_{1n}  \\
%%%%%%%%%%%%%%%%%%%%%%%%%%%%%%%%%%%%%%%%%%%%%%%%%%%%%%%%%%%%%%%%%%%%%%%%%%%%%%
\lambda_{12} + \lambda_{23} & \lambda_{23} & 0 & \lambda_{34} & \cdots & 
\lambda_{34} + \cdots + \lambda_{(n-2)(n-1)} & \lambda_{12} + \lambda_{23} + \lambda_{1n} \\
%%%%%%%%%%%%%%%%%%%%%%%%%%%%%%%%%%%%%%%%%%%%%%%%%%%%%%%%%%%%%%%%%%%%%%%%%%%%%%%
\lambda_{12} + \lambda_{23} + \lambda_{34} & \lambda_{23} + \lambda_{34} &
\lambda_{34} & 0 & \ddots & \vdots & \vdots \\
%%%%%%%%%%%%%%%%%%%%%%%%%%%%%%%%%%%%%%%%%%%%%%%%%%%%%%%%%%%%%%%%%%%%%%%%%%%%%%%
  * &   *   &  *  &  * & \ddots  & \lambda_{(n-2)(n-1)} & \lambda_{(n-2)(n-1)} 
+ \lambda_{(n-1)n} \\ 
%%%%%%%%%%%%%%%%%%%%%%%%%%%%%%%%%%%%%%%%%%%%%%%%%%%%%%%%%%%%%%%%%%%%%%%%%%%%%%
 *   &   *  &  * & * & \cdots &   0   & \lambda_{(n-1)n} \\
%%%%%%%%%%%%%%%%%%%%%%%%%%%%%%%%%%%%%%%%%%%%%%%%%%%%%%%%%%%%%%%%%%%%%%%%%%%%%%
 *   &   *  &  * & * & \cdots &  \lambda_{(n-1)n} &  0   
\end{array}
\right)
$ } 
\end{center}

\begin{center}
\rotatebox[origin=c]{90}{
\tiny
$
\Lambda^4 \ = \
\left(
\renewcommand{\arraystretch}{3}
\begin{array}{cccccccc}
0 & \lambda_{12} & \lambda_{12} + \lambda_{23} & 
\renewcommand{\arraystretch}{1}
\begin{array}{c}
\lambda_{12} + \lambda_{23} \\ 
\mbox{} + \lambda_{34} 
\end{array} & 
\renewcommand{\arraystretch}{1}
\begin{array}{c}
\lambda_{12} + \lambda_{23} \\
\mbox{} + \lambda_{34} + \lambda_{45} 
\end{array} & 
\cdots & \lambda_{12} + \cdots + \lambda_{(n-2)(n-1)} & \lambda_{1n} \\
%%%%%%%%%%%%%%%%%%%%%%%%%%%%%%%%%%%%%%%%%%%%%%%%%%%%%%%%%%%%%%%%%%%%%%%%%%%%%%
\lambda_{12} & 0 & \lambda_{23} & \lambda_{23} + \lambda_{34} & 
\renewcommand{\arraystretch}{1}
\begin{array}{c}
\lambda_{23} + \lambda_{34}  \\
\mbox{} + \lambda_{45} 
\end{array} & 
\cdots & \lambda_{23} + \cdots + \lambda_{(n-2)(n-1)} &  \lambda_{12} + \lambda_{1n}  \\
%%%%%%%%%%%%%%%%%%%%%%%%%%%%%%%%%%%%%%%%%%%%%%%%%%%%%%%%%%%%%%%%%%%%%%%%%%%%%%
\lambda_{12} + \lambda_{23} & \lambda_{23} & 0 & \lambda_{34} & 
\lambda_{34} + \lambda_{45} & \cdots & \lambda_{34} + \cdots + \lambda_{(n-2)(n-1)} & 
\renewcommand{\arraystretch}{1}
\begin{array}{c}
\lambda_{12} + \lambda_{23}  \\
\mbox{} + \lambda_{1n} 
\end{array} \\
%%%%%%%%%%%%%%%%%%%%%%%%%%%%%%%%%%%%%%%%%%%%%%%%%%%%%%%%%%%%%%%%%%%%%%%%%%%%%%%
\renewcommand{\arraystretch}{1}
\begin{array}{c}
\lambda_{12} + \lambda_{23} \\ 
\mbox{} + \lambda_{34} 
\end{array} & 
\lambda_{23} + \lambda_{34} &
\lambda_{34} & 0 & \lambda_{45} & \cdots & \lambda_{45} + \cdots + \lambda_{(n-2)(n-1)} & 
\renewcommand{\arraystretch}{1}
\begin{array}{c}
\lambda_{12} + \lambda_{23}  \\
\mbox{} + \lambda_{34} + \lambda_{1n}
\end{array}  \\
%%%%%%%%%%%%%%%%%%%%%%%%%%%%%%%%%%%%%%%%%%%%%%%%%%%%%%%%%%%%%%%%%%%%%%%%%%%%%%
\renewcommand{\arraystretch}{1}
\begin{array}{c}
\lambda_{12} + \lambda_{23}  \\
\mbox{} + \lambda_{34} + \lambda_{45} 
\end{array} & 
\renewcommand{\arraystretch}{1}
\begin{array}{c}
\lambda_{23} + \lambda_{34} \\
\mbox{} + \lambda_{45} 
\end{array} & 
\lambda_{34} + \lambda_{45}& \lambda_{45} & 0 & \ddots & \vdots & \vdots \\
%%%%%%%%%%%%%%%%%%%%%%%%%%%%%%%%%%%%%%%%%%%%%%%%%%%%%%%%%%%%%%%%%%%%%%%%%%%%%%%
  * &   *   &  *  &  * & * & \ddots  & \lambda_{(n-2)(n-1)} & \lambda_{(n-2)(n-1)} 
+ \lambda_{(n-1)n} \\ 
%%%%%%%%%%%%%%%%%%%%%%%%%%%%%%%%%%%%%%%%%%%%%%%%%%%%%%%%%%%%%%%%%%%%%%%%%%%%%%
 *   &   *  &  * & * & * & \cdots &   0   & \lambda_{(n-1)n} \\
%%%%%%%%%%%%%%%%%%%%%%%%%%%%%%%%%%%%%%%%%%%%%%%%%%%%%%%%%%%%%%%%%%%%%%%%%%%%%%
 *   &   *  &  * & * & * & \cdots &  \lambda_{(n-1)n} &  0   
\end{array}
\right)
$ } 
\end{center}

\begin{center}
\rotatebox[origin=c]{90}{
\tiny
$
\Lambda^k \ = \
\left(
\renewcommand{\arraystretch}{3}
\begin{array}{ccccccccc}
0 & \lambda_{12} & \lambda_{12} + \lambda_{23} & 
\renewcommand{\arraystretch}{1}
\begin{array}{c}
\lambda_{12} + \lambda_{23} \\ 
\mbox{} + \lambda_{34} 
\end{array} & 
\renewcommand{\arraystretch}{1}
\begin{array}{c}
\lambda_{12} + \lambda_{23} \\
\mbox{} + \lambda_{34} + \lambda_{45} 
\end{array} & 
\cdots & \cdots & 
\renewcommand{\arraystretch}{1}
\begin{array}{c}
\lambda_{12} + \cdots  \\
\mbox{} + \lambda_{(n-2)(n-1)} 
\end{array} & \lambda_{1n} \\
%%%%%%%%%%%%%%%%%%%%%%%%%%%%%%%%%%%%%%%%%%%%%%%%%%%%%%%%%%%%%%%%%%%%%%%%%%%%%%
\lambda_{12} & 0 & \lambda_{23} & \lambda_{23} + \lambda_{34} & 
\renewcommand{\arraystretch}{1}
\begin{array}{c}
\lambda_{23} + \lambda_{34}  \\
\mbox{} + \lambda_{45} 
\end{array} & 
\cdots & \cdots & 
\renewcommand{\arraystretch}{1}
\begin{array}{c}
\lambda_{23} + \cdots \\ 
\mbox{} + \lambda_{(n-2)(n-1)} 
\end{array} &  \lambda_{12} + \lambda_{1n}  \\
%%%%%%%%%%%%%%%%%%%%%%%%%%%%%%%%%%%%%%%%%%%%%%%%%%%%%%%%%%%%%%%%%%%%%%%%%%%%%%
\lambda_{12} + \lambda_{23} & \lambda_{23} & 0 & \lambda_{34} & 
\lambda_{34} + \lambda_{45} & \cdots & \cdots &
\renewcommand{\arraystretch}{1}
\begin{array}{c}
\lambda_{34} + \cdots \\
\mbox{} + \lambda_{(n-2)(n-1)} 
\end{array} & 
\renewcommand{\arraystretch}{1}
\begin{array}{c}
\lambda_{12} + \lambda_{23}  \\
\mbox{} + \lambda_{1n} 
\end{array} \\
%%%%%%%%%%%%%%%%%%%%%%%%%%%%%%%%%%%%%%%%%%%%%%%%%%%%%%%%%%%%%%%%%%%%%%%%%%%%%%%
\vdots & \vdots & \vdots & \ddots & \ddots & \cdots & \cdots &\vdots&\vdots \\
%%%%%%%%%%%%%%%%%%%%%%%%%%%%%%%%%%%%%%%%%%%%%%%%%%%%%%%%%%%%%%%%%%%%%%%%%%%%%%
\renewcommand{\arraystretch}{1}
\begin{array}{c}
\lambda_{12} + \lambda_{23} + \\ 
\mbox{} \cdots  + \lambda_{(k-1)k} 
\end{array} &
\renewcommand{\arraystretch}{1}
\begin{array}{c}
\lambda_{23} + \lambda_{34} +  \\
\mbox{} \cdots + \lambda{(k-1)k} 
\end{array} &
\cdots & \lambda_{(k-1)k} & 0 & \lambda_{k(k+1)}& \cdots & 
\renewcommand{\arraystretch}{1}
\begin{array}{c}
\lambda_{k(k+1)} + \cdots \\
+ \lambda_{(n-2)(n-1)}
\end{array} & 
\renewcommand{\arraystretch}{1}
\begin{array}{c}
\lambda_{12} + \lambda_{23} + \cdots  \\
\mbox{} + \lambda_{(k-1)k} + \lambda_{1n}
\end{array}  \\
%%%%%%%%%%%%%%%%%%%%%%%%%%%%%%%%%%%%%%%%%%%%%%%%%%%%%%%%%%%%%%%%%%%%%%%%%%%%%%
\vdots & \vdots & \vdots & \vdots & \vdots & \ddots&\ddots & \vdots & \vdots \\
%%%%%%%%%%%%%%%%%%%%%%%%%%%%%%%%%%%%%%%%%%%%%%%%%%%%%%%%%%%%%%%%%%%%%%%%%%%%%%
  * &   *   &  *  &  * & * & \cdots & \ddots & \lambda_{(n-2)(n-1)} &
\renewcommand{\arraystretch}{1}
\begin{array}{c}
 \lambda_{(n-2)(n-1)} + \\
\mbox{} + \lambda_{(n-1)n} 
\end{array} \\ 
%%%%%%%%%%%%%%%%%%%%%%%%%%%%%%%%%%%%%%%%%%%%%%%%%%%%%%%%%%%%%%%%%%%%%%%%%%%%%%
 *   &   *  &  * & * & * & \cdots & \cdots &  0   & \lambda_{(n-1)n} \\
%%%%%%%%%%%%%%%%%%%%%%%%%%%%%%%%%%%%%%%%%%%%%%%%%%%%%%%%%%%%%%%%%%%%%%%%%%%%%%
 *   &   *  &  * & * & * & \cdots & \cdots & \lambda_{(n-1)n} &  0   
\end{array}
\right)
$ } 
\end{center}

\newpage

%%%%%%%%%%%%%%%%%%%%%%%%%%%%%%%%%%%%%%%%%%%%%%%%%%%%%%%%%%%%%%%%%%%%%%%

\vspace{1cm}

\noindent
Escuela de Matem\'aticas \newline
Universidad Industrial de Santander \newline
Apartado Aereo 678, Bucaramanga \newline
Colombia \newline
{\sc Email}: mparedes@uis.edu.co


\begin{thebibliography}{99}

\bibitem{Be} {\sc A. Besse}, {\it Einstein Manifolds}, Ergb. Math. 
Grenzgeb. (3) 10, Springer Verlag, Berlin, Heilderberg, New York, 1987. 

\bibitem{Bl} {\sc M. Black}, {\it Harmonic Maps into Homogeneous Spaces}, 
Pitman Res. Notes Math. Ser., vol. 255, Longman, Harlow, 1991.

\bibitem{Bo} {\sc A. Borel}, {\it K\"ahlerian Coset Spaces of 
Semi--Simple Lie Groups}, Proc. Nat. Acad. of Sci. {\bf 40} (1954), 
1147-1151.

\bibitem{BH} {\sc A. Borel \& F. Hirzebruch}, {\it Characteristic 
classes and homogeneous spaces I}, Amer. J. Math. {\bf 80} (1958), 
458-538.

\bibitem{BS} {\sc F. E. Burstall \& S. Salamon}, {\it Tournaments, 
Flags and Harmonic Maps}, Math. Ann. {\bf 277} (1987), 249-265.

\bibitem{ChW} {\sc S. S. Chern \& J. G. Wolfson}, {\it Harmonic Maps 
of the Two--Sphere into a Complex Grassmann Manifold II}, Ann. of Math. 
{\bf 125} (1987), 301-335.

\bibitem{ES} {\sc J. Eells \& J. H. Sampson}, {\it Harmonic Mappings 
of Riemannian Manifolds}, Amer. J. Math. {\bf 86} (1964), 109-160.

\bibitem{ESa} {\sc J. Eells \& S. Salamon}, {\it Twistorial 
Constructions of Harmonic Maps of Surfaces into Four--Manifolds}, 
Ann. Scuola Norm. Sup. Pisa (4) {\bf 12} (1985), 589-640.

\bibitem{EW} {\sc J. Eells \& J. C. Wood}, {\it Harmonic Maps from 
Surfaces to Complex Projective Spaces}, Advances in Mathematics {\bf 49} 
(1983), 217-263.

\bibitem{GH} {\sc A. Gray \& L. M. Hervella}, {\it The Sixteen 
Classes of Almost Hermitian Manifolds and Their Linear Invariants}, 
Ann. Mat. Pura Appl. {\bf 123} (1980), 35-58.

\bibitem{L} {\sc A. Lichnerowicz}, {\it Applications Harmoniques et 
Va\-ri\'e\-t\'es K\"ahl\'e\-rie\-nnes}, Symposia Mathematica 3 (Bologna,
1970), 
 341-402.

\bibitem{M} {\sc J. W. Moon}, {\it Topics on Tournaments}, Holt, 
Rinehart and Winston, New York, 1968.

\bibitem{MN1} {\sc X. Mo \& C. J. C. Negreiros}, 
{\it (1,2)--Symplectic Structures on Flag Manifolds}, Tohoku Mathematical
Journal, V. 57, No. 02 (2000), 271-283.  

 
\bibitem{MN2} {\sc X. Mo \& C. J. C. Negreiros}, {\it Tournaments 
and Geometry of Full Flag Manifolds}, Proceedings of the XI Brazilian 
Topology Meeting, Rio Claro, Brazil, World Scientific (2000).

\bibitem{N1} {\sc C. J. C. Negreiros}, {\it Some Remarks about Harmonic 
Maps into Flag Manifolds}, Indiana University Mathematics Journal 
{\bf 37}, No. 3 (1988).

\bibitem{N2} {\sc C. J. C. Negreiros}, {\it Harmonic Maps from Compact 
Riemann Surfaces into Flag Manifolds}, Thesis, University of Chicago, 
1987.

\bibitem{P1} {\sc M. Paredes}, {\it Aspectos da Geometria Complexa
das Variedades Bandeira}, Doctoral Thesis, Universidade Estadual de
Campinas, Brasil, 2000.

\bibitem{P2} {\sc M. Paredes}, {\it 
Some Results on the Geometry of Full 
Flag Manifolds and Harmonic Maps}, pre--print.

\bibitem{Sa} {\sc S. Salamon}, {\it Harmonic and Holomorphic Maps}, 
Lecture Notes in Mathematics 1164, Springer, 1986.

\bibitem{WG} {\sc J. A. Wolf \& A. Gray}, {\it Homogeneous Spaces 
Defined by Lie Groups Automorphisms}. II, Journal of Differential 
Geometry, Vol. 2, No. 2, (1968), 115 - 159.

\end{thebibliography}
\end{document}